\documentclass[12pt]{article}
\usepackage{amsmath,amssymb,amsfonts,latexsym}
\usepackage[dvips]{epsfig}
\usepackage{amscd}

\title{Link homology and Frobenius extensions}
\author{Mikhail Khovanov} 
\date{November 23, 2004}

\newtheorem{prop}{Proposition}

\newtheorem{corollary}{Corollary}

\begin{document}
\maketitle
\baselineskip 14pt

\begin{abstract} We explain how rank two Frobenius extensions 
of commutative rings lead to link homology theories and discuss 
relations between these theories, Bar-Natan theories, 
equivariant cohomology and the Rasmussen invariant.  
\end{abstract} 

\def\C{\mathbb C}
\def\Z{\mathbb Z}
\def\Q{\mathbb Q}
\def\F{\mathbb F}
\def\S{\mathbb S}
\def\l{\lbrace}
\def\r{\rbrace}
\def\o{\otimes}
\def\lra{\longrightarrow}
\newcommand{\dmod}{\mathrm{-mod}}
\def\cF{\mathcal{F}}
\def\Hom{\textrm{Hom}}
\def\drawing#1{\begin{center} \epsfig{file=#1} \end{center}}

\begin{center} \emph{AMS Subject Classification: } 57M27 \end{center} 

\vspace{0.15in}

{\bf Frobenius systems.} 
Suppose $\iota: R \lra A$ is an inclusion of commutative rings, and $\iota(1)=1.$ 
The restriction functor 
$$\mbox{Res}: A\dmod \lra R\dmod$$ 
has left and right adjoint functors: the induction functor 
$\mbox{Ind}(M)= A\otimes_R M$ and
the coinduction functor $\mbox{CoInd}(M) = \Hom_R(A,M).$ 
Following Kadison [Ka] and others, we say that $\iota$ is a Frobenius extension if 
the induction and coinduction functors are isomorphic. Equivalently, $\iota$ 
is Frobenius if the restriction functor has a biadjoint (two-sided adjoint). 
We note that Kadison [Ka] treats the more general case 
of not necessarily commutative $R$ and $A.$ In this paper we consider 
only commutative rings. The following proposition 
is well-known (or see [Ka, Section 4]). 

\begin{prop} $\iota$ is a Frobenius extension iff there exists an $A$-bimodule 
map $\Delta : A \lra A\otimes_R A$ and an $R$-module map $\epsilon: A \lra R$ 
such that $\Delta$ is coassociative and cocommutative, and 
 $(\epsilon \otimes \mathrm{Id})\Delta = \mathrm{Id}.$
\end{prop} 

A Frobenius extension, together with a choice of $\epsilon$ and $\Delta,$ 
will be denoted $\cF= (R,A,\epsilon, \Delta)$ and called a 
 \emph{Frobenius system} (as in [Ka]). 

A Frobenius system defines a 2-dimensional TQFT, a tensor functor 
from oriented $(1+1)$-cobordisms to $R$-modules, by assigning $R$ to 
the empty 1-manifold, $A$ to the circle, $A\otimes_R A$ to the 
disjoint union of two circles, etc. The structure maps $\iota, \epsilon, 
\Delta, m$ (where $m$ is the multiplication in $A$) are assigned to 
basic two-dimensional cobordisms, see [A], [Ka], or [Kh1]. 

Each numbered formula (1), (2), etc. in this paper describes a Frobenius 
system. We denote the Frobenius system associated with formula 
($i$) by $\cF_i=(R_i, A_i, \epsilon, \Delta).$

\vspace{0.1in} 

{\bf Link homology.} In [Kh1] we constructed an invariant of links, based on a particular 
Frobenius system that we denote here by $\cF_1,$ with $R_1= \Z,$  
$A_1= \Z[X]/(X^2),$ the counit 
and comultiplication 
\begin{equation} \label{justZ} 
\epsilon(1)=0, \hspace{0.1in} \epsilon(X)=1, \hspace{0.1in} 
\Delta(1) = 1\otimes X + X \otimes 1, \hspace{0.1in} \Delta(X) = X\otimes X. 
\end{equation}  
To a plane diagram $D$ of a link $L$ we 
assigned a commutative cube of $R_1$-modules and $R_1$-module 
homomorphisms, then passed to the total complex of the cube and 
took its cohomology, which ended up being independent 
of the choice of diagram $D.$ The construction worked as well for 
Frobenius system $\cF_2$ with $R_2=\Z[c],$  $A_2=\Z[X,c]/(X^2)$ and  
\begin{equation} \label{Zandc} 
 \begin{array}{lll} 
   \epsilon(1)=-c,  &  & \Delta(1) = 1\otimes X + X\otimes 1 + c X \otimes X, \\
   \epsilon(X)= 1, & & \Delta(X) = X \otimes X. 
 \end{array} 
\end{equation} 
This system is graded, with $\deg(X)=2, \deg(c)=-2$ (we follow grading 
conventions of [Kh2] rather than those of the earlier paper [Kh1]).  
The multiplication, comultiplication, 
unit and counit maps  have degrees $0,2 ,0$ and $-2$ respectively. As a result, 
homology theories associated with this Frobenius system and with $\cF_1$ are
 bigraded. 

\vspace{0.1in} 

In a recent paper [BN] Dror Bar-Natan defined several new homology theories 
of links and tangles. 
One of his insights, when viewed from an algebraic viewpoint, is in 
modifying the equation $X^2=0$ to $X^2=t,$ where $t$ is a formal 
variable (equal to one-eighth of his invariant of a closed genus 3 
surface). Namely, a certain quotient of his link invariant [BN, Section 9.2] is 
the homology theory assigned to Frobenius system $\cF_3$ 
with $R_3=\Z[t],$   $A_3=\Z[X],$  $X^2=t$ (i.e., $\iota: R_3 \to A_3$ 
maps $t$ to $X^2$) and 
 \begin{equation} \label{Zandt} 
\begin{array}{lll} 
 \epsilon(1)=0 , & &  \Delta(1)=1\otimes X + X\otimes 1,   \\
 \epsilon(X)=1 , & &  \Delta(X)= X\otimes X + t 1\otimes 1. 
\end{array}  
 \end{equation} 
This theory is graded, with $\mathrm{deg}(t) =4.$ The invariant of 
a link is a complex of graded free $\Z[t]$-modules, up to chain homotopy, 
and its homology is a bigraded link homology theory. 

\vspace{0.1in} 

{\bf Base change.} 
Let $\cF$ be a Frobenius system. Any homomorphism  $\psi:R \to R'$
of commutative rings with $\psi(1)=1$ determines a Frobenius system
$\cF'=(R', A', \epsilon', \Delta')$ where $A'=A\otimes_R R'$ and the 
comultiplication, counit maps are induced from those maps for $A$ by 
tensoring with the identity 
endomorphism of $R'.$  We will say that $\cF'$ is obtained by base change 
from $\cF$ via $\psi.$ For instance, Frobenius system $\cF_1$ 
is obtained by base change (with surjective $\psi$) from each of 
$\cF_2,$  $\cF_3.$ Specifically, adding the relation  $t=0$ 
to $\cF_3$ recovers the $\cF_1$ theory (here $\psi: \Z[t]\to \Z,$ 
$\psi(t)=0$), 
 while the ring homomorphism $\Z[t]\lra \Q$ taking $t$ to $1$ 
produces Frobenius extension $\Q\subset \Q[X]/(X^2-1)$ 
leading to Lee's theory [L] and an application to slice genus [R]. 
The specialization $t=1$ collapses the grading.

\vspace{0.1in} 

{\bf The dual system} of $\cF$ is the Frobenius system 
$\cF^{\ast}=(R,A^{\ast},\iota^{\ast},m^{\ast}),$ 
where $A^{\ast}= \Hom_R (A,R)$ and the unit, multiplication, counit, 
comultiplication 
of $\cF^{\ast}$ are obtained by dualizing the counit, comultiplication, unit, 
multiplication of $(R,A),$ respectively.  We say that $\cF$ is self-dual if 
$\cF^{\ast}$ is isomorphic to $\cF$ as a Frobenius system, i.e. 
there is a $R$-module isomorphism $A^{\ast}\cong A$ which 
identifies $\iota^{\ast}$ with $\epsilon,$ $m^{\ast}$ with $\Delta,$ etc.

\begin{prop} Frobenius systems $(R_i, A_i)$ for $i=1,2,3$ are self-dual. 
\end{prop} 
This is proved by a direct computation. Note that if a system is self-dual, so are  
its base changes. 

\vspace{0.1in} 

{\bf Twisting.} Given a Frobenius system $\cF$ and an invertible 
element $y\in A,$ we can twist the comultiplication and counit by $y$:
$$ \epsilon'(x) = \epsilon(yx), \hspace{0.1in} \Delta'(x) = \Delta(y^{-1}x).$$ 
This results in a Frobenius system $\cF'=(R,A, \epsilon', \Delta').$ 
Twisting by invertible elements of $A$ is the only way to modify 
the counit and comultiplication in Frobenius extensions, see [Ka, Theorem 1.6].

\vspace{0.1in} 

To any Frobenius system $\cF$ and a link diagram $D$ 
we can assign a complex of $R$-modules, denoted $\cF(D),$ using the
algorithm of [Kh1]. 

\begin{prop} Complexes $\cF(D)$ and $\cF'(D)$ are isomorphic if 
$\cF'$ is a twisting of $\cF.$ 
\end{prop}
\emph{Proof:} $\cF(D)$ is the total complex of an $n$-cube $V_{D,\cF},$
assigned to $D$ and $\cF,$ where $n$ is the number of crossings of $D.$ 
In each vertex of the cube stands a tensor power of $A,$ and each arrow 
is either a multiplication or a comultiplication map on some factors 
of $A^{\otimes k}.$ 

In fact, $V_{D,\cF}$ and $V_{D, \cF'}$ are 
isomorphic, as cubes of $R$-modules. The isomorphism is constructed by 
making it the identity on the source vertex of the cube and extending it 
arrow by arrow to the whole cube. For each vertex, the isomorphism 
$A^{\otimes k}  \lra A^{\otimes k}$ is a multiple of the identity map, 
with the coefficient being a product of powers of $y_i^{-1},$ $1\le i\le k,$ where 
$$y_i = 1\otimes \dots \otimes 1\otimes y \otimes 1\otimes \dots \otimes 1\in 
A^{\otimes k}.$$    
This isomorphism of cubes induces an isomorphism of complexes $\cF(D)$ 
and $\cF'(D).$     $\square$ 

\begin{corollary} If $\cF'$ is a twisting of $\cF$ then $\mathrm{H}(D,\cF) \cong 
\mathrm{H}(D, \cF').$ 
\end{corollary} 

Let's look at the $\cF_2$ theory from this viewpoint. 
The element $1+cX\in A_2$ is invertible, $(1+cX)(1-cX)=0.$ 
The twisting $\cF_2'$ of $\cF_2$ by this element is a Frobenius 
system with the counit and comultiplication given by formula (\ref{justZ}). 
Therefore, $\cF_2'$ is just the Frobenius system $\cF_1$ with the 
ground ring extended by adding a free variable $c.$

\begin{corollary} For any link $L$ and its diagram $D$ complexes 
$\cF_1(D) \otimes_{\Z} \Z[c]$ and $\cF_2(D)$ are isomorphic. Consequently,  
there is an isomorphism of cohomology groups 
 $$ \mathrm{H}(L, \cF_2) \cong \mathrm{H}(L, \cF_1) \otimes_{\Z} \Z[c].$$ 
\end{corollary} 

Therefore, adding $c$ does not provide us with new information.
Moreover,  Jacobsson [J] showed that the $\cF_2$ theory is not invariant 
under cobordisms. Multiplication by $1+cX$ appears in his work as well. 

Thus, twisting preserves link homology groups but not, in general, link 
homology functors. 

\vspace{0.1in} 

{\bf Invariance under the first Reidemeister move.} 
We would like to know for which Frobenius systems $\cF$  
the complexes $\cF(D_1)$ and $\cF(D_2)$ are chain homotopy equivalent 
whenever $D_1,D_2$ are related by a Reidemeister move. Any such system 
produces a homology theory of links (where to a link $L$ we assign homology
 $\mathrm{H}(D,\cF)$ of the complex $\cF(D)$ for a diagram $D$ of $L$).   

The complex assigned to a one-crossing knot diagram is either 
$$ 0 \lra A\otimes_R A \stackrel{m}{\lra}  A \lra 0 $$ 
or 
$$ 0 \lra A \stackrel{\Delta}{\lra} A\otimes_R A \lra 0 $$
depending on the crossing's sign. Each of these complexes must be 
chain homotopy equivalent to $0 \lra A \lra 0,$ via an $R$-linear, and, 
preferably, $A$-linear homotopy. The unit map $\iota\otimes \mathrm{Id}$ 
is a section of $m,$ and the first complex decomposes as a sum of a contractible 
complex and $0 \lra \mathrm{ker}(m) \lra 0.$ 
Hence, we need an $A$-module isomorphism 
$\mathrm{ker}(m) \cong A.$ Since 
$$\mathrm{ker}(m)\cong \mathrm{coker}(\iota) \otimes_R A
\cong (A/R)\otimes_R A$$ 
(where $A/R$ is the $R$-module quotient), the invariance follows if  $A/R$ is 
a free $R$-module of rank $1.$ Assuming this and pulling back $1\in R \cong A/R$ to 
$A$ we get $X\in A$ such that $A\cong R1 \oplus RX$ as $R$-module.  

We say that $\cF$ has rank two if there exists $X\in A$ such that 
$A\cong R1 \oplus RX.$ 

\vspace{0.1in} 

{\bf A universal rank two Frobenius system.}  Let $\cF_4$ be the Frobenius system with 
$$R_4 = \Z[a,c,e,f,h,t]/(ae-cf, af+chf-cet-1)$$
and $A_4=R_4[X]/(X^2-hX-t).$ The comultiplication and counit are
\begin{equation} \label{universal} 
\begin{array}{l} 
 \Delta(1)= (et-hf)1\otimes 1 + e X\otimes X + 
  f(1\otimes X + X\otimes 1), \\
 \Delta(X)= ft 1\otimes 1 + et (1\otimes X + X\otimes 1) + 
 (f+eh) X\otimes X, \\
   \epsilon(1)= -c, \hspace{0.2in} \epsilon(X)= a. 
 \end{array} 
\end{equation} 
 $R_4$ and $A_4$ are graded, with $a,c,e,f,h,t,X$ in degrees 
$0,-2,$  $-2,0,$  $2,$ $4,2$ respectively. 

\begin{prop} \label{any} $\cF_4$ is a rank two graded 
Frobenius system, universal 
in the following sense. Suppose $\cF'$ is a rank two Frobenius system, and 
$X'\in A'$ a splitting element, $A'\cong R'1\oplus R'X'.$ There exists a unique 
homomorphism $\psi: A_4\lra A'$ with $\psi(X)=X'$ that realizes $\cF'$ as 
a base change of $\cF_4.$ 
\end{prop} 

\emph{Proof:}  That $\cF_4$ is a graded rank two Frobenius system can be 
verified by a direct computation. Furthermore, for $X'$ as above, we have 
$X'^2 = h' X'+ t', \epsilon'(1)=-c'$ and $\epsilon'(X')=a'$ for 
unique $h',t',a',c'\in R'.$ Write 
$$ \Delta'(1) = d' 1\otimes 1 + e' X'\otimes X' + f'(1\otimes X' + X'\otimes 1).$$
Since $\Delta'$ is a map of $A'$-bimodules, 
$$ (X' \otimes 1 ) \Delta'(1) = (1\otimes X')\Delta'(1)$$ 
implying $d'=e' t'-h'f'.$ Define $\psi$ by taking generators $a,c,e,f,h,t,X$ of $A_4$ to 
the corresponding elements of $A'.$   $\square$ 
 
\vspace{0.1in} 

\emph{Remark:} 
This system has a geometric description in terms of dotted surfaces (where a 
dot stands for multiplication by $X$). For instance, the invariant of a sphere 
decorated by a dot is $\epsilon(X)=a.$ The relation $X^2=hX+t$ 
translates into the following skein relation: a disk decorated by two 
dots equals a disk decorated by a single dot times $h$ plus  
a disk times $t.$ The formula for $\Delta(1)$ translates into 
a surgery skein relation for a tube (compare with [Kh3] and [BN, Section 11.2]). 
Adding a dot to a surface corresponds to taking the connected sum 
of the surface with the torus and dividing by $2,$ in the language of 
Bar-Natan [BN].  

\vspace{0.1in} 

Since $f+eX$  is invertible, with the inverse $a+ch-cX,$ 
 and has degree $0,$ we can twist by it. The comultiplication and counit 
become 
\begin{equation}  \label{univ}
  \begin{array}{lll} 
   \epsilon(1)=0,  &  &  \Delta(1) =1\otimes X + X\otimes 1 - h 1\otimes 1    \\
   \epsilon(X)=1,  &  &  \Delta(X)= X\otimes X + t \hspace{0.04in} 1\otimes 1. 
  \end{array}
\end{equation} 
After the twist, all the structure maps depend on $h$ and $t$ only. 
Let $\cF_5$ be the Frobenius system with 
$$ R_5 = \Z[h,t], \hspace{0.15in} A_5 = R_5[X]/(X^2-hX-t), \hspace{0.15in} 
 \deg(h)=2, \hspace{0.1in} \deg(t)=4,$$
and $\epsilon, \Delta$ given by (\ref{univ}).  
We have 
\begin{prop} \label{anyrank2} Any rank two Frobenius system is obtained from 
$\cF_5$ by a composition of base change and twist. 
\end{prop} 

There is a map from the Bar-Natan's theory [BN] to $\cF_5$ 
given by introducing dotted surfaces and skein relations 
between them. Coefficients of  skein relations are the  
structure constants of $\cF_5.$ The decomposition 
of $\Delta(1)$ in formula (\ref{univ}) becomes a skein 
relation for dotted surfaces that implies the $4Tu$ relation of Bar-Natan. 
Although this map from Bar-Natan's theory to $\cF_5$ is 
neither surjective nor injective already for the empty link, 
his arguments are universal, being instantly adoptable to $\cF_5$ 
and any base change of the latter. The next result follows from 
[BN] at once. 

\begin{prop} Complexes $\cF_5(D_1)$ and $\cF_5(D_2)$ are isomorphic 
if $D_1, D_2$ are two diagrams of the same oriented link. Frobenius 
system $\cF_5$ determines a bigraded link homology theory which is 
functorial for link cobordisms up to sign indeterminancy. Any (graded) 
Frobenius system obtained by (graded) base change from $\cF_5$ 
determines a (bi)graded homology theory of links, functorial for 
link cobordisms up to sign indeterminancy.  
\end{prop} 

In view of proposition~\ref{anyrank2}, any rank two Frobenius 
system (after twisting, if necessary) produces a cohomology theory 
of links. On the chain level, all the information is already captured by 
$\cF_5,$ since $\cF'(D) \cong \cF_5(D) \otimes_R R'$ for 
any $\cF'$ given by base change from $\cF_5.$ 

Bar-Natan's theory [BN] cannot be immediately 
expressed in the language of Frobenius systems.  
Indeed, for Frobenius systems the homology of 2-component 
unlink is $A\otimes_R A,$ the second tensor power of the unknot's homology $A$ over 
the homology $R$ of the empty knot. In Bar-Natan's framework there is no such 
isomorphism. In particular, a tube (which should represent $\Delta(1)$) cannot 
be cut into a sum of surfaces separating the two boundaries of the tube. 

\vspace{0.1in} 

{\bf Examples.} 

{\bf 1.}  The quotient of $\cF_5$ by the ideal $(h)$ produces Frobenius system 
$\cF_3$ described earlier. Variable $h$ can also be removed by $X\to X-\frac{h}{2}$ 
if $2$ is made invertible in the ground ring (at the cost of modifying $t$).   

\vspace{0.07in} 

{\bf 2.} Frobenius system $\cF_5$ is almost self-dual. Namely, the dual of 
$\cF_5$ is the system where $h$ is changed to $-h$ in the structure maps. 
In particular, any base change $\psi$ of $\cF_5$ with $\psi(2h)=0$  is  
self-dual. Dual Frobenius systems make an appearance in knot homology, for 
there is an isomorphism of chain complexes  
 $$ C(D^!, \cF) \cong C(D, \cF^{\ast})^{\ast} $$ 
where $D^!$ is the mirror image of $D.$  
If $\cF$ is a rank two Frobenius extension, its dual also has rank two, and 
the above equation descends to an isomorphism of link invariants. 

\vspace{0.07in} 

{\bf 3.} Link homology theory discovered by Bar-Natan in [BN, Section 9.3] 
and investigated by Turner [T] is given by a Frobenius system $\cF_6$ with 
$$ R_6= \mathbb{F}_2[H], \hspace{0.1in} A_6= R_6[X]/(X^2-HX), \hspace{0.1in} \deg(H)=2,$$
where $\mathbb{F}_2$ is the 2-element field, and
 \begin{equation} \label{over2} 
 \begin{array}{lll} 
   \epsilon(1) = 0,  & & \Delta(1) = 1 \otimes X + X\otimes 1 + H 1\otimes 1, \\
   \epsilon(X)= 1, & &  \Delta(X) = X\otimes X. 
 \end{array} 
 \end{equation} 
This system is self-dual. The base change $\psi: R_5\lra R_6$ to this 
system is given by $ \psi(h)=H,$ $\psi(t)=0.$ 

\vspace{0.1in}

{\bf 4. Field extensions.} Any field extension $R\subset A$ of finite degree 
is Frobenius, and any nonzero $R$-linear map $A\lra R$ can 
serve as a counit $\epsilon.$ In particular, any degree two field extension 
$R\subset A$ gives rise to a link homology theory. We do the base change 
$R\lra \overline{R},$ where $\overline{R}$ is the algebraic closure of $R.$ 
This base change preserves the dimension of homology groups 
(as $R$- , respectively, $\overline{R}$-vector spaces): 
$$ \mathrm{dim}_R \hspace{0.05in}\mathrm{H}(L,\cF) = 
\mathrm{dim}_{\overline{R}} \hspace{0.05in} \mathrm{H}(L, \overline{\cF}),$$
where $\overline{\cF}$ is the above base change of $\cF=(R,A,\epsilon, \Delta).$    
Let $\overline{A}=A\otimes_R \overline{R}.$
There are two cases to consider. 

\vspace{0.07in} 

{\bf a.} The extension $R\subset A$ is separable. Then 
$\overline{A}\cong \overline{R} \times \overline{R},$ as an $\overline{R}$-algebra, 
and the resulting theory is the one studied by Lee [L] in characteristic $0,$ 
and by Shumakovitch [Sh] in finite characteristic, with 
$$\mathrm{dim}_{\overline{R}} \hspace{0.05in}\mathrm{H}
(L, \overline{\cF})= 2^{m},$$
where $m$ is the number of components of $L.$ In particular, the total 
rank of homology groups $\mathrm{H}(L,\cF)$ depends only on $m.$ 

\vspace{0.07in} 

{\bf b.} The extension $R\subset A$ is inseparable. Then $\mathrm{char}R=2$ 
and $A\cong R[y]/(y^2+t)$ where $t\in R, \sqrt{t} \notin R.$ We have 
$\overline{A} \cong \overline{R}[X]/(X^2)$ where $X=y + \sqrt{t},$ and 
$\sqrt{t}\in \overline{R}.$ Twisting $\overline{\epsilon},$ we can assume 
$\overline{\epsilon}(1)=0,$ $\overline{\epsilon}(X)=1.$ 
The resulting theory $\overline{\cF}$ is obtained from $\cF_1$ by base change 
$$ \Z \lra \mathbb{F}_2 \lra \overline{R}.$$ 
The intermediate theory (with $R= \mathbb{F}_2$) is simply 
the original theory of [Kh1, Section 7] with coefficients in the 
2-element field. Denote it by $\mathrm{H}(L, \mathbb{F}_2).$ We have 
 $$\mathrm{dim}_{R}\hspace{0.05in} \mathrm{H}(L, \cF) = 
 \mathrm{dim}_{\mathbb{F}_2}\hspace{0.05in} \mathrm{H}(L,\mathbb{F}_2)$$ 
for any $\cF$ given by an inseparable degree two field extension $R\subset A.$ 

\vspace{0.1in} 

\begin{figure} [htb] \drawing{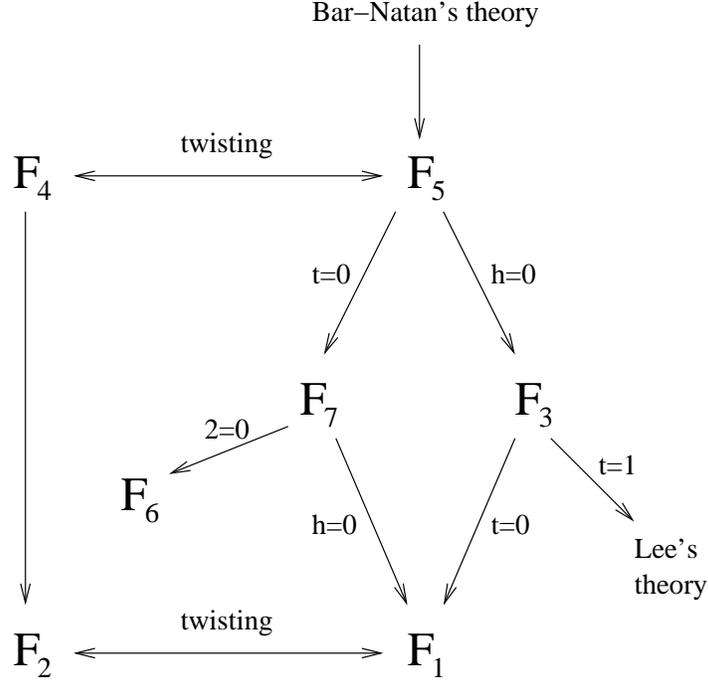}\caption{A diagram of link homology 
theories} \label{first-fig} 
 \end{figure}

\vspace{0.1in} 

{\bf Equivariant cohomology.} Frobenius extension $\cF_1$ has a cohomological 
interpretation: $R_1\cong \Z$ is isomorphic to the cohomology ring of a 
point, and $A_1\cong \Z[X]/(X^2)$  to the cohomology ring 
of a 2-sphere. The trace map $\epsilon$ is the integration along the
fundamental cycle on $\mathbb{S}^2.$ 

Other extensions that we considered have similar interpretations via 
equivariant cohomology. Suppose that a topological group $G$ acts continuously 
on $\mathbb{S}^2.$  Define 
$$ R_G \cong \mathrm{H}^{\ast}_G(p, \Z) = \mathrm{H}^{\ast}(BG,\Z)$$ 
to be the $G$-equivariant cohomology ring of a point $p$ (where $BG$ is 
the classifying space of $G$), and 
$$ A_G \cong \mathrm{H}^{\ast}_G(\mathbb{S}^2, \Z) = 
 \mathrm{H}^{\ast}(\mathbb{S}^2\times_G EG,\Z)$$ 
the equivariant cohomology ring of the 2-sphere. 
Then, in several cases, $(R_G, A_G)$ is a rank $2$ Frobenius extension, 
with $\epsilon$ induced by the integration along the fibers of the 
$\mathbb{S}^2$-fibration 
 $$  \mathbb{S}^2\times_G EG \lra BG.$$ 

{\bf Examples.}  

{\bf 1.} The standard action of $G=SU(2)$ on $\mathbb{C}^2$ induces 
an action on $\mathbb{S}^2$ (with $-I$ acting trivially) 
and leads to the Frobenius system $\cF_3.$ Indeed, 
 \begin{eqnarray*} 
 R_3  & = & \Z[t] \cong \mathrm{H}^{\ast}(BSU(2), \Z) =  
 \mathrm{H}^{\ast}(\mathbb{HP}^{\infty}, \Z),  \\
 A_3 & = & \Z[X] \cong \mathrm{H}^{\ast}(\mathbb{S}^2 \times_{SU(2)} 
   ESU(2), \Z)\cong \mathrm{H}^{\ast}(\mathbb{CP}^{\infty},\Z), \hspace{0.2in} 
  X^2=t. 
 \end{eqnarray*} 
$X$ here is the two-dimensional cohomology class of $\mathbb{CP}^{\infty}$ 
which evaluates to $1$ on $\mathbb{CP}^1\subset \mathbb{CP}^{\infty}.$ 
We choose $t\in  \mathrm{H}^{4}(\mathbb{HP}^{\infty}, \Z) $
so that its pullback to $\mathbb{CP}^{\infty}$ equals $X^2.$ Also,  $\epsilon(X)=1.$ 
The geometric counterpart of the $SU(2)$-equivariant theory was 
considered by Seidel and Smith [S]. 

\vspace{0.1in} 

{\bf 2.} Taking  $G$ to be the group $U(2)$ with the usual action on $\mathbb{S}^2$
(so that the center $U(1)$ acts trivially), we get Frobenius system 
$\cF_5$: 
\begin{eqnarray*}
   R_5  & = & \Z[h,t] \cong  \mathrm{H}^{\ast}(BU(2), \Z) \cong 
 \mathrm{H}(\mathrm{Gr}(2, \infty), \Z),  \\ 
   A_5 & = & \Z[h,X] \cong \mathrm{H}^{\ast}(\mathbb{S}^2 \times_{U(2)} 
   EU(2), \Z) \cong\mathrm{H}^{\ast}(BU(1)\times BU(1), \Z)  . 
 \end{eqnarray*} 
$\mathrm{Gr}(2, \infty)$ is the Grassmannian of complex planes $\mathbb{C}^2$ 
in $\mathbb{C}^{\infty},$ its cohomology ring freely generated by $h$ and $t$ 
in degrees $2$ and $4,$ while $BU(1)\cong \mathbb{CP}^{\infty}.$ 
   Notice that $A_5$ is the polynomial ring with generators 
$X,Y=h-X,$ and $R_5$ is the ring of symmetric functions in $X$ and $Y,$ 
with $h$ and $-t$ the elementary symmetric functions. 

\vspace{0.1in} 

{\bf 3.} $G=U(1),$ the group of rotations of $\mathbb{S}^2$ about a 
fixed axis. In this case we get a Frobenius system $\cF_7$ which is the 
quotient of $\cF_5$ by the ideal $(t)$: 
 \begin{equation}
  \begin{array}{lll}  
  R_7  =  \Z[h],  & &  A_7 = \Z[X,h]/(X^2-hX),   \\
  \epsilon(1) = 0,  &   & \Delta(1) = 1 \otimes X + X\otimes 1 - h 1\otimes 1, \\
   \epsilon(X) = 1,  &   & \Delta(X) = X \otimes X, \\
   \end{array} 
  \end{equation} 
The system $\cF_6$ is the modulo $2$ specialization of $\cF_7$ and can be described 
via equivariant $U(1)$-cohomology with coefficients in $\mathbb{F}_2.$  

If we change the coefficient ring of equivariant cohomology from $\Z$ to 
any field $k$ of characteristic other than $2,$ then substitution $X\to X-\frac{h}{2}$ 
takes us to the theory which is a  base change of $\cF_3$ (the latter consider 
over $k$ too) with $t=-\frac{h^2}{4}.$  The homology of the resulting theory 
is given by suitably doubling (with a shift) the $\cF_3$ homology 
(over $k$).

\begin{figure} [htb] \drawing{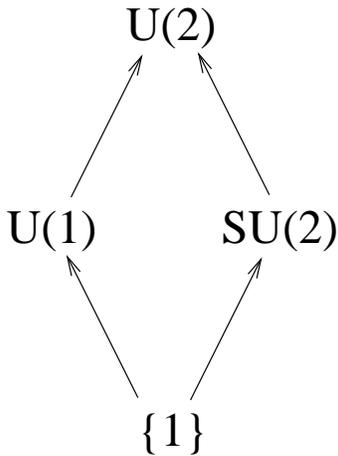}\label{diag}\caption{Some 
  subgroups of $U(2)$}
 \end{figure}

Each of the four theories $\cF_1,\cF_3, \cF_5, \cF_7$ in the vertices of the 
central rombus in figure~\ref{first-fig} is the $G$-equivariant theory for some 
connected closed subgroup of $U(2),$ see figure~2. Group 
inclusion arrows in figure~2 are reversals of base change arrows 
in figure~\ref{first-fig}. 

\vspace{0.1in} 

{\bf Bar-Natan's theory and the Rasmussen invariant.} To a diagram $D$ of an 
oriented link $L$ we can associate the complex 
$\cF_3(D)$ of graded free $\Z[t]$-modules (recall that $t$ is 
one-eighth of the closed genus three surface in  [BN], as well as 
a generator of $\mathrm{H}^4(\mathbb{HP}^{\infty},\Z)$). 
 The chain 
homotopy equivalence class of this complex is an invariant 
of $L.$  Denote the complex 
$\cF_3(D)\otimes_{\Z} \Q$ of graded free $\Q[t]$-modules 
by $C_t(D)$ and its cohomology by $\mathrm{H}_t(D)$ (this is 
exactly Bar-Natan's universal theory over $\Q$).  
The cohomology $ \mathrm{H}_t(D)$ 
is a finitely-generated graded $\Q[t]$-module and thus a direct 
sum of torsion modules $\Q[t]/(t^m),$ for various $m,$ and
free modules $\Q[t].$ Denote by $\mathrm{Tor}(D)$ the 
torsion submodule of $\mathrm{H}_t(D)$ and by $\mathrm{H}'(D)$ the 
quotient $\mathrm{H}_t(D)/\mathrm{Tor}(D).$ We denote corresponding 
link homologies by $\mathrm{Tor}(L)$ and $\mathrm{H}'(L),$ 
respectively. Both of these theories are 
cobordism-friendly: given a link cobordism 
$S$ from $L_1$ to $L_2$ there are well-defined (up to 
overall sign) homomorphisms  
 \begin{eqnarray*}  
    & &  \mathrm{Tor}(S) : \mathrm{Tor}(L_1) \lra 
    \mathrm{Tor}(L_2), \\
    & &      \mathrm{H}'(S)  :  \mathrm{H}'(L_1) \lra 
     \mathrm{H}'(L_2). 
 \end{eqnarray*}
 To define $\mathrm{H}'(S)$ on an element
 $\alpha\in \mathrm{H}'(L_1),$ 
pull $\alpha$ back to $\mathrm{H}(L_1)$ and  
apply the homomorphism $\mathrm{H}_t(S)$ 
composed with the quotient map $\mathrm{H}_t(L_2) 
\to \mathrm{H}'(L_2).$ 
The map $\mathrm{Tor}(S)$ is simply the restriction 
of $\mathrm{H}_t(S)$ to the torsion submodule of $\mathrm{H}(L_1).$  

\begin{prop} $\mathrm{H}'(L)$ is a free $\Q[t]$-module of rank 
$2^m,$ where $m$ is the number of components of $L.$ 
\end{prop}

\emph{Proof.} This follows Lee [L]. The quotient of the complex 
$C_t(D)$ by its subcomplex $(t-1)C_t(D)$ has cohomology 
of dimension equal to the rank of $\mathrm{H}'(L).$  
$\square$ 
 
\vspace{0.1in} 

Assume now $L$ is a knot. Then $C_t(D)$ is naturally 
a complex of free $\Q[X]$-modules, where $X^2=t.$ 

\begin{prop} $\mathrm{H}'(L)$ is a free $\Q[X]$-module of
 rank one concentrated in cohomological degree $0.$ 
\end{prop} 

Therefore, as a graded $\Q[X]$-module, 
$\mathrm{H}'(L) \cong \Q[X] \{ -s'(L)-1 \}$ for some even integer 
$s'(L).$ Notice that $\mathrm{H}_t(L)$ is nontrivial only 
in odd $q$-degrees, since $L$ is a knot, hence the 
shift is by an odd degree. It is clear from the definition 
that $s'(L)$ is the Rasmussen invariant $s(L),$ assuming that 
we normalize as in [Kh2,4].   

The quotient of $\mathrm{H}'(L)$ by the ideal $(t-1)$ gives 
Lee's theory. In the latter any connected 
knot cobordism induces a nontrivial homomorphism, see [R].  
Lifting to $\mathrm{H}'(L)$, we have 

\begin{prop} Any connected genus $g$ cobordism $S$ between 
knots $L_1$ and $L_2$ induces a nonzero grading-preserving 
map $\mathrm{H}'(S): \mathrm{H}'(L_1) \to \mathrm{H}'(L_2)\{-2g\}.$ 
\end{prop} 

All other results of Rasmussen [R] admit a natural interpretation 
via $H'$ as well. Frobenius system $\cF_3$ allows working with graded 
rather than filtered complexes, making Rasmussen structures 
slightly more explicit and bundling up a number of invariance 
results from several papers into that for $\cF_3$ implied 
by Bar-Natan [BN]. 

\vspace{0.1in} 

The complex $C_t(D)/tC_t(D)$ is the original complex of [Kh1, Section 7] 
with coefficients in $\Q,$ which we denote $C(D).$  Since 
$\Q[X]$ has homological dimension $1,$ the complex $C_t(D)$ 
is isomorphic to the direct sum of complexes 
 $$ 0 \lra \Q[X]\{2m+i\} \stackrel{X^m}{\lra} \Q[X]\{i\} \lra 0 \hspace{1in} 
 (\ast )$$
for various $m>0,i\in \Z,$ a contractible complex, and the 
complex 
$$ 0 \lra \Q[X]\{-s(L)-1\} \lra 0.$$ 
  Reducing modulo $t,$ we see that  $C(D)$ is isomorphic, modulo contractible 
complexes, to the direct sum of  
 $$ 0 \lra A\{2m+i\} \stackrel{0}{\lra} A\{i\} \lra 0,  $$  
 for $m>1,$ 
 $$ 0 \lra A\{2+i\} \stackrel{X}{\lra} A\{i\} \lra 0  $$  
($m=1$ case), and the complex $0 \lra A\{-s(L)-1\} \lra 0.$ 
Here $A=\Q[X]/(X^2).$ 

In particular, complexes 
$$ 0 \lra A \stackrel{X}{\lra} A \stackrel{X}{\lra} 
 \dots \stackrel{X}{\lra} A \lra 0 $$ 
of length greater than $1$ (i.e. with more than two $A$'s) 
cannot appear as direct summands of $C(D),$ 
which settles one of the problems implicitly raised in [Kh4, Section 3]. 
Thus, the rank of reduced homology $\widetilde{\mathrm{H}}(L)$ 
is always less than the rank of $\mathrm{H}(L),$ 
with the difference 
$\dim \mathrm{H}(L) - \dim \widetilde{\mathrm{H}}(L)-1$ being twice 
the number of terms in $C_t(D)$ of the form $(\ast)$ for $m>1.$

\vspace{0.1in} 

{\bf Acknowledgments.} I am grateful to Dror Bar-Natan and Dmitry Fuchs for 
enlightening discussions. Many useful suggestions from the referee 
were incorporated into the second version of the paper. 
While writing this paper, I was partially supported by 
the NSF grant DMS-0407784. 

\vspace{0.1in} 

{\bf References} 

\vspace{0.1in} 

\noindent 
[A] L.~Abrams, Two dimensional topological quantum field theories and 
Frobenius algebras, \emph{J. Knot Theory and its Ramifications} 5 (1996) 569-587. 

\noindent 
[BN] D.~Bar-Natan, Khovanov's homology for tangles and cobordisms, 
arXiv math.GT/0410495. 

\noindent 
[J] M.~Jacobsson, Khovanov's conjecture over Z[c], math.GT/0308151. 

\noindent 
[Ka] L.~Kadison, New examples of Frobenius extensions, University Lecture 
Series  14, AMS, 1999. 

\noindent
[Kh1] M.~Khovanov, A categorification of the Jones polynomial,  
 \emph{Duke Math. J.} 101 (2000), no. 3, 359--426, math.QA/9908171.

\noindent 
[Kh2] M.~Khovanov, A functor-valued invariant of tangles,   
\emph{Algebr. Geom. Topol.} 2 (2002), 665--741 (electronic), math.QA/0103190.

\noindent
[Kh3] M.~Khovanov, sl(3) link homology I, 
\emph{Algebr. Geom. Topol.} 4 (2004) 1045-1081, math.QA/0304375. 

\noindent 
[Kh4] M.~Khovanov, Patterns in knot cohomology I, \emph{Experimental 
Math.} 12 (2003), no.~3, 365--374, math.QA/0201306

\noindent 
[L] E.~S.~Lee, An endomorphism of the Khovanov invariant, 
to appear in \emph{Adv. Math.,} arxiv 
math.GT/0210213. 

\noindent 
[R] J.~Rasmussen, Khovanov homology and the slice genus,  
math.GT/0402131. 

\noindent 
[Sh] A.~Shumakovitch, Torsion of the Khovanov homology, 
math.GT/0405474. 

\noindent 
[S] I.~Smith, Talk at MSRI, March 2004. 

\noindent 
[T] P.~Turner, Calculating Bar-Natan's characteristic two Khovanov homology, 
math.GT/0411225. 

\vspace{0.2in} 

\noindent
Department of mathematics, University of California, Davis CA 95616

\noindent 
mikhail@math.ucdavis.edu

\end{document}